\documentclass[amsfonts,12pt]{article}

\usepackage{amsfonts}
\usepackage{amsmath}
\usepackage{amssymb}
\usepackage{array, rotating}

\newtheorem{theorem}{Theorem}
\newtheorem{corollary}[theorem]{Corollary}

\newtheorem{definition}[theorem]{Definition}
\newtheorem{example}[theorem]{Example}

\newtheorem{remark}{Remark}
\def\ppp{{\mathbb{P}}}

\def\zzz{\mathbb{Z}}

\def\pf{{\bf proof}:\ }
\def\qed{$\Box$}
\DeclareMathOperator{\Aut}{Aut}

\DeclareMathOperator{\Perm}{Perm}

\begin{document}

\author{David Joyner and Amy Ksir\thanks{Mathematics Dept, USNA,
Annapolis, MD 21402,
wdj@usna.edu and ksir@usna.edu}
}
\title{Automorphism groups of some AG codes}
\date{1-19-2005}

\maketitle

\begin{abstract}

We show that in many cases, the automorphism group of a curve and
the permutation automorphism group of a corresponding AG code are
the same.  This generalizes a result of Wesemeyer \cite{W} beyond
the case of planar curves.
\end{abstract}

\vskip .5in

\section{Introduction}

The construction of AG codes uses the Riemann-Roch space $L(D)$
associated to a divisor $D$ of a curve $X$ defined over a finite
field \cite{G}. Typically $X$ has no non-trivial automorphisms, but
when it does we may ask how this can be used to better understand AG
codes constructed from $X$.

Conversely, we may ask how the permutation automorphism group of an
AG code corresponds with the automorphism group of the curve used to
construct the code. In this paper we show that, in many cases, the
automorphism group of a curve and the permutation automorphism group
of a corresponding AG code are in fact the same.

Knowledge of which codes have large automorphism group can have
applications to encoding (see \cite{HLS}) and to decoding (indeed,
permutation decoding has been implemented in version 2.0 of the
error-correcting computer algebra package \cite{GUAVA}).

\section{The Riemann-Roch space $L(D)$ and the associated AG code.}
\label{sec:1}

Let $X$ be a smooth projective curve (scheme of dimension 1) over a
finite field $F$, and let $F(X)$ denote the field of rational
functions on $X$. If $D$ is any divisor on $X$, the Riemann-Roch
space $L(D)$ is a finite dimensional $F$-vector space given by
\[
L(D)=L_X(D)= \{f\in F(X)^\times \ |\ div(f)+D\geq 0\}\cup \{0\},
\]
where $div(f)$ denotes the (principal) divisor of the function $f\in
F(X)$. If $\bar{D}$ denotes the corresponding divisor over the
algebraic closure $\bar{F}$, then $L(\bar{D}) = L(D) \otimes
\bar{F}$ \cite{Sti}, \cite{TV}.

Let $P_1, \ldots ,P_n\in X(F)$ be distinct points, and let
$E=P_1+\ldots +P_n$ be the associated divisor. Let $D$ be a divisor
of positive degree on $X$ such that $D$ and $E$ have disjoint
support. Let $C=C(D,E)$ denote the AG code

\begin{equation}
\label{eqn:AGcode} C=\{(f(P_1),\ldots ,f(P_n))\ |\ f\in L(D)\}.
\end{equation}
This is the image of $L(D)$ under the evaluation map

\begin{equation}
\label{eqn:eval}
\begin{array}{c}
eval_E:L(D)\rightarrow F^n,\\
f \longmapsto (f(P_1),\ldots ,f(P_n)).
\end{array}
\end{equation}

The kernel of the map $eval_E$ is contained in $L(D-E)$, which is
empty if $n > deg(D)$.  Thus for $n > deg(D)$, $eval_E$ defines an
isomorphism between $L(D)$ and the code $C(D,E)$.

\section{From curve automorphisms to code automorphisms.}
\label{sec:codes}

Now let $G$ be a group of automorphisms of the curve $X$, and assume
that $D$ and $E$ are both stabilized by $G$.  We will say that $G
\subseteq \Aut_{D,E}(X)$. Then $G$ also acts on the code $C$, as
follows.

The action of ${\rm Aut}(X)$ on $F(X)$ is defined as:
\[
\begin{array}{cccc}
&{\rm Aut}(X)&\longrightarrow &{\rm Aut}(F(X)),\\
 & T &\longmapsto & (f\longmapsto T^{*}f)
\end{array}
\]
where for any $P \in X$, $T^{*}f(P)=f(T^{-1}(P))$.  (We use $T^{-1}$
rather than $T$ here, to conform to the convention that the action
should be on the left.)

Note that $Y=X/G$ is also smooth and $F(X)^G=F(Y)$.

Of course, ${\rm Aut}(X)$ also acts on the group $Div(X)$ of
divisors of $X$, denoted $T(\sum_P d_P P)=\sum_P d_P T(P)$, for
$T\in \Aut(X)$, $P$ a prime divisor, and $d_P \in \zzz$. It is easy
to see that $div(T^{*}f)=T^{-1}(div(f))$. Because of this, if
$div(f)+D\geq 0$ then $div(T^{*}f)+T^{-1}(D)\geq 0$, for all $T\in
{\rm Aut}(X)$. In particular, if the action of $G\subset {\rm
Aut}(X)$ on $X$ leaves $D\in Div(X)$ stable then $G$ also acts on
$L(D)$. Assuming that $n
> \deg D$, the isomorphism $eval_E : L(D) \to C$ will send this
action to an action of $G$ on $C$.  Specifically, each $T \in G$
acts by

\begin{eqnarray*}
(f(P_1), f(P_2), \ldots, f(P_n)) & \mapsto & (T^{*}f(P_1),
T^{*}f(P_2), \ldots ,T^{*}f(P_n)) \\
& = & (f(T^{-1}(P_1)), f(T^{-1}(P_2)), \ldots, f(T^{-1}(P_n))).
\end{eqnarray*}

 If we also assume that $G$ leaves $E$
stable, then $G$ acts by permutations on the set $\{P_1, \ldots, P_n
\}$. Thus $(T^{-1}(P_1), T^{-1}(P_2), \ldots ,T^{-1}(P_n))$ is a
permutation of the points $(P_1, P_2, \ldots ,P_n)$, and the above
action on $C$ simply permutes the corresponding coordinates.

\begin{definition}
The \textbf{permutation automorphism group} $\Perm C$ of the code
$C\subset F^n$ is the subgroup of $S_n$ (acting on $F^n$ by
coordinate permutation) which preserves $C$.
\end{definition}

Thus if $n > \deg D$, we have defined a homomorphism from
$\Aut_{D,E}(X)$ to $\Perm C$:
\begin{eqnarray*}
\rho: \Aut_{D,E}(X) & \to & \Perm C. \\
T & \mapsto & eval_E \circ T^{*} \circ eval_E^{-1}
\end{eqnarray*}
In the next section, we will construct an inverse for this
homomorphism.

\section {From code automorphisms to curve automorphisms.}

Now we would like to answer the question: when does a group of
permutation automorphisms of the code $C$ induce a group of
automorphisms of the curve $X$?  We will show that permutation
automorphisms of the code $C(D,E)$ induce curve automorphisms
whenever $D$ is very ample and the degree of $E$ is large enough.
Under these conditions, the groups $\Aut_{D,E}(X)$ and $\Perm C$ are
isomorphic.  In proving these facts, we generalize some results of
Wesemeyer \cite{W}, who dealt with the planar case.

\begin{theorem}
\label{thrm:lift} 
Let $X$ be an algebraic curve, $D$ be a very ample
divisor on $X$, and $P_1 \ldots P_n$ be a set of points on $X$
disjoint from the support of $D$.  Let $E = P_1 + \ldots + P_n$ be
the associated divisor, and $C = C(D,E)$ the associated AG code. Let
$G$ be the group of permutation automorphisms of $C$.  Then there is
an integer $r \geq 1$ such that if $n > r \cdot \deg(D)$, then $G$
can be lifted to a group of automorphisms of the curve $X$ itself.
This lifting defines a group homomorphism $\psi: \Perm C \to
\Aut(X)$.  Furthermore, the lifted automorphisms will preserve $D$
and $E$, so the image of $\psi$ will be contained in
$\Aut_{D,E}(X)$.
\end{theorem}

\pf First, note that since $n > \deg D$, $eval_E: L(D) \to C$ is a
vector space isomorphism.  Thus the permutation action of $G$ on $C$
can be pulled back to a linear action on $L(D)$.  Next, we use $D$
to embed $X$ into projective space $\ppp^{d}$, where $d=\dim
L(D)-1$.  If we let $Y_0, \ldots, Y_d$ be a basis for $L(D)$, then
the embedding is given explicitly by

\[
\phi:X\rightarrow \ppp^d,
\]
\[
P\longmapsto [Y_0(P):\ \ldots \ :Y_{d}(P)].
\]
The vector space action of $G$ on $L(D)$ induces an action on the
polynomial ring $F[Y_0, \ldots ,Y_d]$ and a projective linear action
on $\ppp^d$. We will show that under the stated hypotheses, this
action preserves the image of $X$ in $\ppp^d$, so restricts to an
action on $X$.  Furthermore, this action will stabilize the divisors
$D$ and $E$.

To prove these claims, let us look more carefully at the action of
$G$. Let $\tau$ be an element of $G$; it acts by a permutation of
the coordinates of a point in $C$.  The pullback of $\tau$ to $L(D)$
is the composition $eval_E^{-1} \circ \tau \circ eval_E$, which
by abuse of notation we denote again by $\tau$.  In the
middle, $\tau$ acts as

\begin{equation}
\label{eqn:1}
(f(P_1), \ldots, f(P_n)) \mapsto (f(P_{\tau(1)}), \ldots,f(P_{\tau(n)}))
\end{equation}
where $f$ was a function in $L(D)$.  Because the permutation action
leaves the code $C$ invariant, this new point is also in the code.
Since $eval_E$ is an isomorphism, there is a function which we will
call $\tau(f)$ in $L(D)$ such that

\begin{equation}
\label{eqn:*1}
(f(P_{\tau(1)}), \ldots, f(P_{\tau(n)}))=(\tau(f)(P_1), \ldots,
\tau(f)(P_n)).
\end{equation}
This defines the action of $G$ on $L(D)$.  In particular, the action
of $G$ on the basis $Y_0, \ldots, Y_d$ of $L(D)$ defines an action
of $G$ on polynomial ring $F[Y_0, \ldots, Y_d]$:

\begin{equation}
\label{eqn:*4}
Y_0^{e_0} \ldots Y_d^{e_d}\mapsto (\tau Y_0)^{e_0} \ldots (\tau Y_d)^{e_d}.
\end{equation}
Then the action on
the projective space $\ppp^d$ is as follows:  an element $\tau$ of $G$
will act on a point $[Y_0: \ldots :Y_d]$ in $\ppp^d$ via

\begin{equation}
\label{eqn:*2}
\tau[Y_0:\  \ldots \ :Y_d] = [\tau^{-1}Y_0:\  \ldots \ :
\tau^{-1}Y_d].
\end{equation}

Now we will consider how this action on the projective space affects
the images of the points $P_1, \ldots, P_n$ under the embedding
$\phi$. For each point $P_i$, its image $\phi(P_i)$ has projective
coordinates $[Y_0(P_i):\ \ldots \ :Y_d(P_i)]$.  Then

\begin{equation}
\label{eqn:*}
\begin{array}{lll}
\tau^{-1}(\phi(P_i)) & = & \tau^{-1}[Y_0(P_i):\ \ldots \ :Y_d(P_i)] \\
& = & [\tau Y_0(P_i):\ \ldots \ :\tau Y_d(P_i)] \\
& = & [Y_0(P_{\tau(i)}):\ \ldots \ :Y_d(P_{\tau(i)})] \\
& = & \phi(P_{\tau(i)}).
\end{array}
\end{equation}
Thus, a permutation of the code acts by the inverse permutation on
the images of the points of $E$.

Now we would like to show that the image $\phi(X)$ is preserved by
the action (\ref{eqn:*2}) of $G$ on $\ppp^d$. In the case where $d=1$, $X$ must
have genus $0$ and $\phi$ is an isomorphism, so this automatically
holds. For $d>1$, the coordinates $Y_i$ must satisfy some
homogeneous polynomial relations defining $\phi(X)$. Let
$R_1(Y_0,\ldots ,Y_{d})=0, \ldots , R_k(Y_0,\ldots ,Y_{d})=0$ denote
a set of polynomials of minimal degree that define the ideal of
$\phi(X)$ in $\ppp^d$, so that its projective coordinate ring is

\[
F[Y_0, \ldots,Y_{d}]/(R_1,\ldots ,R_k).
\]
Since $R_1, \ldots ,R_k$ are polynomials in $Y_0, \ldots , Y_d$, and
$Y_0, \ldots , Y_d$ are in $L(D)$, the $R_1, \ldots ,R_k$ will be in $L(rD)$
for some $r \geq 1$.
In particular if we let $r$ be the largest degree of the $R_i$'s
in the $Y_j$'s, then each $R_i$ will be in $L(rD)$.  (Often this is true for a
smaller $r$, in fact).

Now let $\tau \in G$ be an automorphism of the code, and consider
the image $\tau(\phi(X))$ of $\phi(X)$ under the induced action 
(\ref{eqn:*2}) on $\ppp^d$. The ideal of $\tau(\phi(X))$ is generated by $\tau(R_1),
\ldots ,\tau(R_k)$, by (\ref{eqn:*4}). 
If we can show that these functions are also in
the ideal of $\phi(X)$, for any $\tau \in G$, then the ideals will
be equal and we will have given an action of $G$ on $\phi(X)$, which
we then pull back via $\phi$ to an action on $X$: for $P\in X$,

\begin{equation}
\label{eqn:*3}
\tau(P) = \phi^{-1}(\tau(\phi(P))).
\end{equation}

Let $R_i\in L(rD)$ be one of the minimal degree generators of the
ideal of $\phi(X)$; $R_i$ is of degree at most $r$ in $Y_0, \ldots
Y_d$. Since the action of $G$ is linear on $L(D)$, $\tau(R_i)$ will
also be of degree at most $r$ in $Y_0, \ldots Y_d$, so $\tau(R_i)$
will also be in $L(rD)$.  Since $R_i$ is in the ideal of $\phi(X)$,
$R_i$ vanishes at every point of $X$, including $P_1, \ldots ,P_n$.
Since as we showed above $\tau$ acts as a permutation of the points
$P_i$, $\tau(R_i)$ must also vanish on $P_1, \ldots, P_n$, so
$\tau(R_i)$ vanishes on $E$. This means that $\tau(R_i)$ is in $L(rD
- E)$. But if $n > r\cdot \deg (D)$, then $rD - E$ is a divisor of
degree $< 0$ and $L(rD - E)$ is the trivial vector space, so
$\tau(R_i)$ must vanish identically on $\phi(X)$. Thus $\tau(R_i)$
is in the vanishing ideal associated to $\phi(X)$, for each $R_i$
and for every $\tau \in G$.

We have shown that the action of $G$ on the code gives an action 
(\ref{eqn:*2}) on
$\phi(X)$, which we then pull back via the embedding to an action on
$X$.  At each stage, the action was multiplicative, so we have a
homomorphism $\psi: \Perm C \to \Aut(X)$.  Using (\ref{eqn:*}) it follows
that $E$ is invariant under this action; we now need to show that the
action leaves $D$ invariant. Consider an element $\tau$ of $G$ and
its action on $D$. Because the action (\ref{eqn:*3}) of $\tau$ on $X$ was defined
via an action (\ref{eqn:*1}) on $L(D)$, we know that $\tau$ preserves $L(D)$. But
suppose that $\tau$ did not preserve $D$ itself, so that $\tau(D) =
D'$, $D \neq D'$, but $L(D) = L(D')$. Then there must be a point $P$
in the support of $D$ such that its coefficient, $d_P$, in $D$ is
larger than its coefficient $d_P'$ in $D'$.  Now consider a function
$f \in L(D)$. Because it is also in $L(D')$, we must have $div(f)+D'
\geq 0$. Thus the coefficient of $div(f)$ at $P$ must be at most
$d_P'$. Thus $div(f) + D - (d_P - d_P')P \geq 0$, so in particular
$f$ is in $L(D-P)$. This is true for any $f$ in $L(D)$, so $L(D) =
L(D-P)$. But we assumed that $D$ was very ample; in particular
$L(D)$ separates points, which means that $\dim L(D-P) = \dim L(D) -
1$, a contradiction.  So the action of $G$ on $X$ must preserve $D$.
This means that the image of the homomorphism $\psi: \Perm C \to
\Aut(X)$ is in $\Aut_{D,E}(X)$. \qed

It should be clear from these constructions that $\rho$ and $\phi$,
when they exist, are inverses of each other, making $\Aut_{D,E}(X)$
and $\Perm C$ isomorphic groups.

The result below is actually slightly stronger than the
corresponding result of Wesemeyer (Corollary 4.9 \cite{W}) for
elliptic curves and elliptic codes.

\begin{corollary}
\label{corollary:main} Let $X$ be a smooth projective curve of genus
$g \geq 2$.  Let $D$ be a divisor on $X$ with $\deg D \geq 2g+1$ and
let $E$ be a collection of at least $(1+g)\deg D$ points on $X$
disjoint from the support of $D$. Then the group of permutation
automorphisms of the code $C=C(D,E)$ is isomorphic to the group of
automorphisms of $X$ that fix both $D$ and $E$.
\end{corollary}

\pf Since $\deg D \geq 2g+1$, $D$ is very ample, so we use Theorem
\ref{thrm:lift}; we want to estimate $r$.  Suppose that the image of
the embedding $|D|:X\hookrightarrow \ppp^{d}$ defined over $F$ is
defined by multivariate polynomial relations $R_1=0$, \ldots,
$R_k=0$ over $F$ of minimal degree.  As noted in the proof of
Theorem \ref{thrm:lift}, we can take $r$ to be the maximal degree of
the polynomials $R_1, \ldots, R_k$.  Let $\bar{D}$ be the associated
divisor over the algebraic closure $\bar{F}$.  By ``base-change'',
we see that the image of the associated embedding
$|\bar{D}|:X\hookrightarrow \ppp^{d}$ defined over $\bar{F}$ is
defined by the same multivariate polynomial relations $R_1=0$, ...,
$R_k=0$ over $F$ (and hence over $\bar{F}$).  Note that since $Y_0,
\ldots, Y_d$ form a basis of $L(D)$, they are linearly independent,
so $2 \leq r$ and $X$ cannot be contained in a hyperplane in
$\ppp^d$.

If $d \geq 3$, Gruson, Lazarsfeld, and Peskine \cite{GLP} (since
$\bar{F}$ is algebraically closed) give the maximum degree of the
$R_i's$ as $\deg D + 1 - d$ in most cases, or $\deg D + 2 - d$ if
$X$ has genus zero and its image is smooth and has a $\deg D + 2 -
d$-secant line. In our case, $d = \dim L(D) -1$ and $D$ is
non-special, so $d = \deg D - g$.  Therefore if $g \geq 2$, we will
have $d \geq 3$ and from \cite{GLP}, $r \leq 1+g$.  \qed

There a few special cases to consider that fall outside of Corollary
\ref{corollary:main}. If $X$ is rational, and $d=1$, then the
embedding is an isomorphism and the automorphism groups are the
same.  If $d=2$, then the embedding is as a plane conic, so $r=2$.
For larger $d$, the theorem of \cite{GLP} holds and shows that $r=2$
(and that $X$ always has a 2-secant line, which is not surprising).
In both of these cases, the groups are isomorphic if $\deg E \geq 2
\deg D$.  If $X$ has genus 1 and is embedded smoothly in $\ppp^2$,
it must be as a cubic so $r=\deg D =3$;  the groups will be
isomorphic if $\deg E \geq 3 \deg D=9$.  Again, for larger $d$
\cite{GLP} holds and shows that $r=2$, so the groups are isomorphic
if $\deg E \geq 2 \deg D$.

\begin{remark}
Under the hypotheses of Corollary 3, the length of $C$ is $n=\deg
E$, dimension is $k=\deg D +1 -g$, and minimum distance $d\geq \deg
E -\deg D$ (see for example Corollary II.2.3 \cite{Sti}).
\end{remark}

\begin{example}
Let $F=GF(49)$ and let $X$ denote the curve defined by

\[
y^2=x^7-x.
\]
This has genus $3$. The automorphism group $Aut_F(X)$ is a central
2-fold cover of $PGL_2(F)$: we have a short exact sequence,

\[
1\rightarrow H \rightarrow Aut_F(X)\rightarrow PGL_2(7)\rightarrow
1,
\]
where $H$ denotes the subgroup of $Aut_F(X)$ generated by the
hyperelliptic involution (which happens to also be the center of
$Aut_F(X)$). For details, see \cite{G}, Theorem 1.

Next, we recall some consequences of \S 3.2 in \cite{JT}. There are
$|X(F)|=2\cdot 7^2-7+1=92$ $F$-rational points\footnote{MAGMA
views the curve as embedded in a weighted projective space, with
weights $1$, $4$, and $1$, in which the point at infinity is
nonsingular.}:

\[
X(F)=\{ P_1=[1:0:1], P_2=[0:0:1], \ldots \}.
\]
The automorphism group does not act transitively on $X(F)$ but has
2 orbits: the orbit $C_1$ of $P_1$ and the orbit $C_2=X(F)-C_1$.
We have $|C_1|=7+1=8$ and $|C_2|=2\cdot 7\cdot (7-1)=84$.


Let $D=mP_1$, $E=X(F)-C_1=\{Q_1,\ldots ,Q_{84}\}$, and let

\[
C=C(D,E)=\{(f(Q_1), \ldots ,f(Q_{84}))\ |\ f\in L(D)\}.
\]
This is an $[n,k,d]$ code over $F$, where $n=deg(E)=84$, $k\leq
dim(L(D))$.

Let $G=Stab(P_1,Aut_F(X))$ denote the stabilizer of $P_1$. Since $E$
is an orbit of the full automorphism group, it will also be
stabilized by $G$, so $G=\Aut_{D,E}(X)$.  The group $G$ is a
non-abelian group of order $2\cdot 7\cdot (7-1)=84$.

According to Corollary \ref{corollary:main}, $\Perm C(D,E)$ will be
isomorphic to $G$ if we choose $m$ so that $\deg D$ is at least
$2g+1=7$ and $\deg E=84$ to be at least $(g+1) \deg D=4 \deg D$.
Since $\deg D = m$, this means that $7 \leq m \leq 21$.

Assuming we choose $m>2g-2=4$, the Riemann-Roch theorem implies
$dim(L(D))=m-2$, so $C$ is an $[84,m-2,\geq 84-m]$-code over
$GF(49)$. Since $G$ fixes $D$ and preserves $E$, it acts on $C$ via

\[
g:(f(Q_1), \ldots ,f(Q_{84}))\longmapsto (f(g^{-1}Q_1),
\ldots ,f(g^{-1}Q_{84})),
\]
for $g\in G$.
\end{example}

\begin{remark}
More generally, for $p>3$ and $p\equiv 3\, \pmod 4$,
the curve $X$ defined\footnote{This curve is embedded
into weighted projective space, where $x$ and $z$ have weight $1$
and $y$ has weight $\frac{p+1}{2}$.} 
by $y^2=x^p-x$ over $F=GF(p^2)$ is associated to an
$[n,k,d]$ code $C$ over $F$, where
$n=2p(p-1)$, $k=m-\frac{p-3}{2}$, 
$d\geq 2p^2-2p-m$, provided $m>p$.
This code is the one-point AG code constructed from the 
divisor $D=mP_1$, where $P_1=[1:0:1]$,
and $E$ is the sum of the points in the
orbit $X(GF(p^2))-X(GF(p))$ (see Proposition 3 of \cite{J}).
When $m=p^2$, the parameters of this code beat the 
Gilbert-Varshamov bound \cite{TV}. 
When $p<m<\frac{4p(p-1)}{p+1}$,
using the above corollary, it can be shown that the 
permutation automorphism group $P$ of $C$ is
isomorphic to the stabilizer of $P_1$ in the automorphism group 
of $X$, which is of size $2p(p-1)$. In this case,
$P$ acts on $C$ as a subrepresentation of the regular
representation. It would be interesting to know 
the decomposition of this representation.

In \cite{J},
it is conjectured that $C$ has a permutation decoding algorithm of 
complexity $O(n)$. For a related discussion (for an AG code
constructed from a different curve), see \cite{L}.
\end{remark}


In some interesting cases, there are not enough rational points on
the curve to apply Theorem \ref{thrm:lift}.

\begin{example}
Again, let $X$ denote the genus 3 curve defined by

\[
y^2=x^7-x,
\]
but this time over $F=GF(7)$.  The automorphism group $Aut_F(X)$ is
now a central 2-fold cover of $PSL_2(F)$: we have a short exact
sequence,

\[
1\rightarrow H \rightarrow Aut_F(X)\rightarrow PSL_2(7)\rightarrow
1,
\]
where as before $H$ denotes the subgroup of $Aut_F(X)$ generated by
the hyperelliptic involution (which happens to also be the center of
$Aut_F(X)$). The following transformations are generating elements of
$G$:

\begin{equation}
\label{eqn:3.2}
\begin{array}{cc}
\gamma_1=
\left\{
\begin{array}{c}
x\longmapsto x,\\
y\longmapsto -y,
\end{array}
\right. ,
&
\gamma_2=\gamma_2(a)=
\left\{
\begin{array}{c}
x\longmapsto a^2x,\\
y\longmapsto ay,
\end{array}
\right.\\
\gamma_3=
\left\{
\begin{array}{c}
x\longmapsto x+1,\\
y\longmapsto y,
\end{array}
\right. ,
&
\gamma_4=
\left\{
\begin{array}{c}
x\longmapsto -1/x,\\
y\longmapsto y/x^{4},
\end{array}
\right. 
\end{array}
\end{equation}
where $a\in F^\times$ is a primitive $6-th$ root of unity.

On this curve there are only $8$ $F$-rational points:

\[
X(F)=\{ P_1=[1:0:0], P_2=[0:0:1], P_3=[1:0:1], \ldots ,
P_8=[6:0:1]\}.
\]
Thus it is impossible to choose $D$ and $E$ so that $\deg D \geq 7$
and $E$ consists of at least $4 \deg D$ distinct rational points.
Let us instead choose $D=mP_1$ and $E$ to be all of the other
rational points as before, and compare $\Perm C$ and
$\Aut_{D,E}(X)$.

The automorphism group acts transitively on $X(F)$; as in the
previous example let $G = \Aut_{D,E}(X) = Stab(P_1,Aut_F(X)),$ the
stabilizer of the point at infinity in $X(F)$.  (All of the
stabilizers $Stab(P_i,Aut_F(X))$ are conjugate to each other in
$Aut_F(X)$, $1\leq i\leq 8$). The group $G$ is a non-abelian group
of order $42$ (In fact, the group $G/Z(G)$ is the non-abelian group
of order $21$, where $Z(G)$ denotes the center of $G$.)
Take the automorphisms
$\gamma_1$, $\gamma_2$ with $a=2$ and $\gamma_3$ as generators of $G$. 
If we identify $S=\{P_2,\ldots ,P_8\}$ with $\{1,2,\ldots ,7\}$ then

\[
\gamma_1 \leftrightarrow (2,7)(3,6)(4,5)=g_1,
\]
\[
\gamma_2 \leftrightarrow (2,5,3)(4,6,7)=g_2,
\]
\[
\gamma_3 \leftrightarrow (1,2,\ldots 7)=g_3.
\]



Let $D=5P_1$, $S=C(F)-\{P_1\}$, and let

\[
C(D,E)=\{(f(P_2), \ldots ,f(P_8))\ |\ f\in L(D)\}.
\]
This is a $[7,3,5]$ code over $F$. In fact, $\dim(L(D))=3$, so the
evaluation map $f\longmapsto (f(P_2), \ldots ,f(P_8))$, $f\in L(D)$,
is injective. Since $G$ fixes $D$ and preserves $E$, it acts on $C$
via

\[
g:(f(P_2), \ldots ,f(P_8))\longmapsto (f(g^{-1}P_2),
\ldots ,f(g^{-1}P_8)),
\]
for $g\in G$.

Let $P$ denote the permutation group of this code. It a group of
order $42$. However, it is not isomorphic to $G$! In fact, $P$ has
trivial center. The (permutation) action of $G$ on this code
implies that there is a homomorphism

\[
\rho :G\rightarrow P.
\]
What is the kernel of this map? There are two possibilities:
either a subgroup of order $6$ or a subgroup of order $21$ (this
is obtained using \cite{GAP} by matching possible 
orders of quotients $G/N$ with
possible orders of subgroups of $P$). 
Indeed, the kernel $ker(\phi)=N=\langle g_2,g_3\rangle$ 
is a non-abelian
normal subgroup of $G=\langle g_1,g_2,g_3\rangle$ of order $21$.

\end{example}

\section{Acknowledgements}
We thank Jessica Sidman for the reference in the proof 
of Corollary
\ref{corollary:main} and Will Traves for many 
helpful conversations.

\end{document}